\documentclass[12pt,leqno, numreferences, noid]{kluwer}
\usepackage{amsopn,amsthm,amscd,verbatim}

\renewcommand\copyrightline{Submitted to the
{\it Annals of Global Analysis and Geometry}}

\DeclareMathSymbol{\varSigma}{\mathord}{letters}{"06}
\DeclareMathSymbol{\varPhi}{\mathord}{letters}{"08}

\def\a,{a\kern-.28em\rlap{\lower.80ex\hbox{${}_`$}}\kern.28em\relax}

\theoremstyle{plain}
\newtheorem{prop}{Proposition}
 \newtheorem{thm}{Theorem}
\newtheorem{lemma}{Lemma}

\theoremstyle{definition}
\newtheorem{exmp}{Example}

\renewcommand\acknowledgementsname{Acknowledgments}
\newcommand{\eqref}[1]{{\rm (\ref{#1})}}
\newcommand\half{\mbox{\small{$\frac{1}{2}$}}}
\newcommand\vts{\mkern1mu}
\newcommand\End{{\rm {End}}\vts }
\newcommand\Ext{{\rm {Ext}}\vts }
\newcommand\Hom{{\rm {Hom}}\vts }
\newcommand\Tor{{\rm {Tor}}\vts }
\newcommand\Sq{{\rm {Sq}}\vts }
\newcommand\Spin{{\rm {Spin}}\vts }
\newcommand\Hm{{\rm {H}}\vts }
\newcommand\Pin{{\rm {Pin}}\vts }
\newcommand\Lpin{{\rm {Lpin}}\vts }
\newcommand\SO{{\rm {SO}}\vts }
\newcommand\Ort{{\rm {O}}\vts }
\newcommand\Cl{{\rm {Cl}}\vts}
\newcommand\U{{\rm {U}}\vts}
\newcommand\GL{{\rm {GL}}\vts}

\DeclareMathOperator{\RP}{\mathbb{R}P}
\DeclareMathOperator{\Ad}{Ad}
\DeclareMathOperator{\id}{id}
\DeclareMathSymbol{\geqslant}{\mathrel}{AMSa}{"3E}

\begin{document}
\begin{article}
\begin{opening}

\title{\bf Pin\(^c\) and Lipschitz structures\\ on products of manifolds}

\author{Marcin
\surname{Bobie\'nski}}
\institute{Dept of Mathematical Methods in Physics, Warsaw University \\
Ho\.za 74, 00682 Warszawa, Poland\\ {\it (E-mail: mbobi@mimuw.edu.pl)}}

\author{Andrzej
\surname{Trautman}}
\institute{Institute of Theoretical Physics, Warsaw University\\
Ho\.za 69, 00681 Warszawa, Poland\\
{\it (E-mail: Andrzej.Trautman@fuw.edu.pl)}}

    \begin{abstract}
The topological condition for the existence of a pin\(^c\)
structure on the product of two Riemannian manifolds
is derived and applied to construct examples of manifolds having
the weaker Lipschitz structure, but no pin\(^c\) structure.
An example of a five-dimensional manifold with this property is given;
it is pointed out that there are no manifolds of lower dimension with
this property.
    \end{abstract}

{Primary: 53C27, 57R15; Secondary: 81R25, 83C60.}

\end{opening}

\section{Introduction}

There are two main approaches to complex spinors on Riemannian manifolds
that are not necessarily orientable.
The first assumes the existence of a pin\(^c\) structure, defined
as a reduction of the bundle  of all orthonormal frames of a Riemannian
manifold \(M\) of dimension \(m\) to
the group \(\Pin^c(m)\). The consideration of a pin\(^c\),
rather than pin, structure is motivated, in part, by the need of
physics to describe {\it charged\/} fermions. The insistence on
pin, rather than spin, structures is also influenced
by physics where it is necessary to consider  transformations
of spinor fields under reflections.
 Spin\(^c\) structures appear
 in the theory of the Seiberg--Witten invariants.
 The second approach is based on the notion of
complex {\it spinor bundles}. If \(m\) is even, then these two aproaches are
equivalent: the spinor  bundle is associated with a pin\(^c\)
structure. If \(m\) is odd and \(M\) is not orientable, then a spinor
bundle is associated with a weaker {\it Lipschitz structure\/}, a notion
introduced by Friedrich and Trautman \cite{FT2000}; it is recalled here
in Section \ref{s:prel}.

In every mathematical category  it is of interest to know for
which pairs of  objects a product is defined. For example, there
is a natural, direct product of every pair of Riemannian
manifolds. If both factors have a spin structure (`are spin'),
then their product is also spin. The product of two Riemannian
manifolds has a pin structure if, and only if, both factors are
pin and at least one of them is orientable \cite{CGT}.

In this paper we consider the question of existence of a pin\(^c\)
structure on the product of two closed (i.e. compact without
boundary) Riemannian manifolds.
The main result  consists in the proof, in Section \ref{s:proof},
 of the following
\begin{thm}
Let \(M_1\) and \(M_2\) be closed Riemannian manifolds.
The Riemannian manifold \(M_1\times M_2\) has a pin\(^c\) structure if,
and only if, both \(M_1\) and \(M_2\) have a pin\(^c\) structure
and there holds one of the following conditions:

\noindent {\,\rm (i)} one of the manifolds is orientable,

\noindent {\rm (ii)} the first Stiefel--Whitney classes of both
manifolds have integral lifts, i.e.
there exist elements \(c_i\in\Hm (M_i,\mathbb{Z})\),
\(i=1,2\), such that
\[
w_1(TM_i)\equiv c_i\bmod 2\quad\mbox{for}\quad i=1,2.
\]
\end{thm}
To make the paper self-contained, we recall in Section \ref{s:prel} some results 
on Lipschitz structures proved in \cite{FT2000}.
In Section \ref{s:ex}, Theorem 1 is used to construct  examples of manifolds
with a Lipschitz structure that do not admit a pin\(^c\) structure.
Among them is a  manifold of dimension 5, this being the lowest dimension
for which this phenomenon can occur.
\section{Preliminaries}\label{s:prel}
\subsection{Notation and terminology}
We use the standard notation and terminology, mainly the one of \cite{ThF2000},
\cite{FT2000} and \cite{Lawson89}. For the convenience of the reader, some of the
relevant definitions are presented below.
The ring of integers modulo \(m\), where \(1<m\in\mathbb{N}\), is denoted by
\(\mathbb{Z}/m\); in particular, \(\mathbb{Z}/2\) is the two-element field.
A real, finite-dimensional vector space \(V\) with a positive-definite
quadratic form is referred to as a
{\it Euclidean space\/}; the orthogonal group of its automorphisms
is  \(\Ort(V)\).
The dual of a finite-dimensional vector space \(S\) over \(\mathbb{C}\)
is  \(S^*=\Hom(S,\mathbb{C})\)
and
\(S\times S^*\ni (\varphi,\omega)\mapsto \langle \varphi,\omega
\rangle\in \mathbb{C} \)
 is the evaluation map. The  vector space of all semi-linear
 maps from \(S^*\) to \(\mathbb{C}\) is denoted by \(\bar{S}\);
 there is a semi-linear map \(S\to\bar{S}\), \(\varphi\mapsto\bar{\varphi}\),
 defined by \(\langle \omega,\bar{\varphi}\rangle=\overline
 {\langle\varphi,\omega\rangle}\) for every \(\omega\in S^*\). The vector spaces
 \(\Hom(\bar{S},\mathbb{C})\) and \(\overline{\Hom(S,\mathbb{C})}\)
 can be identified and are denoted by \(\bar{S}^*\).
All manifolds and maps among them are assumed to be smooth;
\(TM\) denotes the total space of the tangent bundle of the
manifold \(M\). The trivial bundle \(M\times \mathbb{R}^k\)
is denoted by \(\theta^k\).
\subsection{Results from topology}
Our main references for algebraic topology are \cite{Greenberg81} and
\cite{MiSt74}. The standard notation for the homology and cohomology
groups is used; the \(i\)th Stiefel--Whitney class of a vector
bundle \(E\to M\)
is denoted by \(w_i(E)\).
\subsubsection{Obstructions to existence of spinor structures}\label{ss:exist}
 Recall that \(w_1(TM)=0\)
is equivalent to the orientability of \(M\).
Depending on whether the group Pin is defined in terms of
a Clifford algebra based on a quadratic form that is positive-
or negative-definite, one has the notion of a pin\(^+\) or
pin\(^-\) structure.
 A Riemannian manifold \(M\) has a pin\(^+\) (resp., pin\(^-\))
 structure if, and only if, \(w_2(TM)=0\) (resp., \(w_2(TM)+w_1(TM)^2=0\)); it
has a pin\(^c\) structure if, and only if, \(w_2(TM)\) is the
mod 2 reduction of an integral cohomology class
\cite{ThF2000,Lawson89}.

\subsubsection{The Universal Coefficient Theorems}\label{ss:UCT}
These theorems relate homology and cohomology
groups with different coefficients. Let
\(R\) be a ring;
there are natural short exact sequences:
\[
0 \to \Hm_n(X,\mathbb{Z})\otimes R \to \Hm_n(X,R)
\to \Tor(\Hm_{n-1}(X,\mathbb{Z}),R) \to 0
\]
and
\[0 \to \Ext(\Hm_{n-1}(X,\mathbb{Z}),R)\to \Hm^n(X,R)
\to \Hom(\Hm_n(X,\mathbb{Z}),R)\to 0. \]

 When \(R\) is a field, then
\[
\Hm^n(X,R) \cong \Hom(\Hm_n(X,R),R).
\]

\subsection{Pin\(^c\) and Lipschitz groups}
\subsubsection{Clifford algebras and spinor representations}\label{ss:srep}
Let \(V\) be a Euclidean,  \(m\)-dimensional  space, and let
 \(\Cl(V)\)  be its Clifford algebra. We
put
\(\Cl^c(V)=\mathbb{C}\otimes\Cl(V).
\)
A volume element \(\eta\in\Cl^c(V)\) associated with \(V\) is
defined as the product of a sequence
of \(m\) pairwise orthogonal vectors
and normalized so that \(\eta^2=1\); if \(\eta\) is a volume element, then
 so is \(-\eta\).
There is a canonical  antiautomorphism \(a\mapsto a^\dag\)
of \(\Cl^c(V)\)  defined as an \(\mathbb{R}\)-linear
isomorphism of the vector structure such that \(a^\dag=\bar{a}\)
for \(a\in\mathbb{C}\subset\Cl^c(V)\), \(v^\dag=v\) for
 \(v\in V\subset\Cl^c(V)\) and \((ab)^\dag=b^\dag a^\dag\)
for every \(a\) and \(b\in\Cl^c(V)\). The antiautomorphism
defines a unitary group
\[
\U(V)=\{a\in\Cl^c(V)\mid a^\dag a=1\}.
\]

Assume from now on  that \(V\) is of odd dimension \(2n-1\),
\(0<n\in\mathbb{N}\), define \(W\) to be the {\it orthogonal\/} 
sum \(V\oplus \mathbb{R}\) and denote by $u\in W$ a unit vector orthogonal to $V$.
Since the dimension \(2n\) of \(W\) is even,
the algebra \(\Cl^c(W)\) is simple and
 there is a complex, \(2^n\)-dimensional vector space \(S\)
  and an isomorphism
\begin{equation}\label{e:Drepr}
\gamma:\Cl^c(W)\to\End S
\end{equation}
of complex algebras with units, i.e., a faithful and irreducible
`Dirac representation' of the Clifford
algebra in \(S\).
 The map \(a\mapsto \overline{\gamma(a^\dag)^*}
\in\End \bar{S}^*\) is also a faithful irreducible representation
of \(\Cl^c(W)\); the simplicity of the algebra  implies that they are equivalent:
there is an isomorphism \(\varPhi:\bar{S}\to S^*\) such that
\(\overline{\gamma(a^\dag)}=\varPhi^{-1}\gamma(a)^*\varPhi\). By rescaling,
\(\varPhi\) can be made to satisfy \(\varPhi^*=\bar{\varPhi}\);
 the Hermitean form
\begin{equation}\label{e:scp}
(\varphi\vts |\psi)=\langle
\varphi,\varPhi\bar{\psi}\rangle,\quad \varphi,\psi\in S,
\end{equation}
is positive
and  invariant with respect to the action of the group \(\U(W)\)
 in \(S\),
 \((\gamma(a)\varphi\vts |\gamma(a)\varphi)=(
\varphi\vts |\varphi)\)
for every \(a\in \U(W)\) and \(\varphi\in S\).
The scalar product \eqref{e:scp} is used to define the adjoint \(f^\dag\)
of \(f\in\End S\) so that \(\gamma(a^\dag)=\gamma(a)^\dag\) for
every \(a\in\Cl^c(W)\).

 The restriction of \eqref{e:Drepr} to \(\Cl^c(V)\subset\Cl^c(W)\)
 is a faithful, but reducible, representation of \(\Cl^c(V)\) in \(S\):
if \(\eta\) is the volume element associated with \(V\), then the subspaces
\begin{equation}\label{e:dec}
S_\pm=\{\varphi\in S\mid \gamma(\eta)\varphi=\pm\varphi\}
\end{equation}
are invariant for this restriction and \(\gamma |\Cl^c(V)=\sigma_+\oplus
\sigma_-\). The `Pauli representations' \(\sigma_+\)
and \(\sigma_-\) of \(\Cl^c(V)\) in \(S_+\) and \(S_-\),
respectively, are irreducible and inequivalent, but not faithful.
For the purposes of this
paper, it is convenient to refer to \(\gamma\) and \(\gamma |\Cl^c(V)\)
as the {\it spinor representations\/} of the algebras in question.

The notation is chosen so that, taking into
account that \(\gamma\) is an isomorphism, one can {\it identify\/}
\(\Cl^c(W)\) with \(\End S\)
 and this is done in the sequel.
\subsubsection{Pin\(^c\) groups}
 The group \(\Pin(W)\subset\Cl(W)\)
is generated by the set of all unit elements of \(W\) and
\(\Pin^c(W)\subset\Cl^c(W)\) is generated by \(\Pin(W)\) and
\(\U_1\). The adjoint representation of \(\Pin^c(W)\) in \(W\) is given by
\(\Ad(a)v=ava^{-1}\). Similar definitions apply to the groups
\(\Pin(V)\) and \(\Pin^c(V)\).

The dimension of \(W\) being even, there is the exact sequence
\begin{equation}\label{e:pin2n}
1\to\U_1\to\Pin^c(W)\stackrel{\Ad}{\longrightarrow}\Ort(W) \to 1.
\end{equation}
The volume elements of the odd-dimensional space \(V\)
are in the kernel of \(\Ad\)
and the image of \(\Pin^c(V)\) by \(\Ad\)
 is the special orthogonal group \(\SO(V)\).
To describe spinors on a nonorientable manifold, one needs a sequence
like \eqref{e:pin2n}, with the full orthogonal group
as the image of \(\Ad\). In odd dimensions, this
can be achieved by either using, instead of \(\Ad\),
 the {\it twisted adjoint representation\/} \cite{ABS} or
extending the group \(\Pin^c(V)\) to a larger {\it Lipschitz group}.
Only the latter approach is suitable when one considers spinor bundles
as the primitive notion.

\subsubsection{The Lipschitz group}
In the notation of Section \ref{ss:srep}, the Lipschitz group, associated
with the odd-dimensional Euclidean
space \(V\subset W\),  is
\[
\Lpin(V)=\{a\in\GL(S)\mid a^\dag a=1\quad\mbox{and}\quad
aVa^{-1}=V\}.
\]
Clearly, \(\Pin^c(V)\) is a subgroup of \(\Lpin(V)\) and
the homomorphism \(\Ad:\Lpin(V)\to\Ort(V)\) is surjective because
\(u\in\Lpin(V)\) and \(\Ad(u)=-\id_V\).
The kernel of \(\Ad\) is the subgroup
\[\{\half (1+\eta)z_++\half (1-\eta)z_-\mid z_+,\vts
z_-\in \U_1\}\cong\U_1\times\U_1
\]
  so that there is an exact sequence
\[
1\to\U_1\times\U_1\to\Lpin(V)\stackrel{\Ad}{\longrightarrow}\Ort(V)\to 1.
\]
\subsection{Spinor bundles and Lipschitz structures}
\subsubsection{Spinor structures}
Let \(M\) be a Riemannian manifold of dimension \(m\); let
\(V\) be a Euclidean space of dimension \(m\) and denote
by \(P\) the \(\Ort(V)\)-bundle of  orthonormal frames on \(M\).
Let us agree to say that \(G\) is  a {\it spinor group\/}
if it contains \(\Spin(V)\) as a subgroup and there is homomorphism
\(\rho:G\to\Ort(V)\) such that \(\rho(\Spin(V))
=\SO(V)\).
A {\it spinor structure\/} of type \((G,\rho)\) on \(M\) is a reduction \(Q\)
of \(P\) to a spinor group \(G\) characterized
by the maps
 \begin{equation}\label{e:CD}
 \begin{CD}
G      @>>>Q\\
 @V\rho VV  @VV {\chi}V\\
 \Ort(V)       @>>>P    @>>>M
 \end{CD}
 \end{equation}
 such that \(\chi(qa)=\chi(q)\rho(a)\) for \(q\in Q\)
 and \(a\in G\).
 \subsubsection{Spinor bundles}
 One defines the complex Clifford bundle
 associated with a Riemannian manifold \(M\) as
 \[
\Cl^c(TM)=\bigcup_{x\in M}\Cl^c(T_x M)\to M.
 \]
 A complex vector bundle \(\varSigma\to M\), with a homomorphism
  \[\tau:\Cl^c(TM)\to\End\varSigma\] of bundles of algebras
  over \(M\), is said to be a
 bundle of Clifford modules. In particular, a {\it spinor bundle\/}
 is a bundle of Clifford modules such that the restriction of \(\tau\)
 to every fibre is a spinor representation in the sense of Section \ref{ss:srep}.

 A spinor bundle \(\varSigma\) on an odd-dimensional manifold is said to be
 {\it decomposable\/} if there is a non-trivial vector bundle
 decomposition \(\varSigma=
 \varSigma_+\oplus\varSigma_-\) such that \(\tau(\Cl^c(TM))\varSigma_\pm
 \subset\varSigma_\pm\). In \cite{FT2000} it is shown that there holds:
 \begin{prop}
A spinor bundle on an odd-dimensional Riemannian manifold \(M\)
is decomposable if, and only if, \(M\) is orientable.
 \end{prop}
 Note that, for every \(x\in M\), there is a decomposition
 \(\varSigma_x=
 \varSigma_{x+}\oplus\varSigma_{x-}\), as in \eqref{e:dec},
  and \(\tau(\Cl^c(T_xM))\varSigma_{x\pm}
 \subset\varSigma_{x\pm}\), but a {\it global\/} decomposition holds
 only in the orientable case.

 Since, by assumption, \(M\) is paracompact (even compact),
 there is a Hermitean
 scalar product on the fibres of \(\varSigma\to M\),
 \[
(.\mid .)_{\!{_\varSigma}}:\varSigma\times_{\!{_M}}\varSigma\to \mathbb{C}.
 \]
 \subsubsection{Lipschitz structures}
 Let \(M\) and \(V\) be a Riemannian manifold and a Euclidean space, both
 of odd dimension \(m\), respectively.
 A {\it Lipschitz structure\/} on
 \(M\) is defined as a spinor structure \eqref{e:CD} of type \((\Lpin(V),\Ad)\).
A spinor bundle \(\varSigma\to M\) with a Hermitean scalar product
 defines, and is associated with,
a Lipschitz structure  such that the fibre \(Q_x\), \(x\in M\),
consists of all isomorphisms \(q:S\to\varSigma_x\) satisfying
\(qVq^{-1}=\tau(T_xM)\) and \((q\varphi |q\varphi)_{\!{_\varSigma}}=
(\varphi |\varphi)\) for every \(\varphi\in S\). The action
of \(\Lpin(V)\) on \(Q\) is by composition of maps.
\begin{prop}\label{p:FT}
An odd-dimensional manifold \(M\) admits a Lipschitz structure if,
and only if, there is a  vector bundle \(\mathbb{R}^2\to E\to M\)
and \(c\in\Hm^2(M,\mathbb{Z})\) such that
\begin{equation}\label{e:FT}
w_2(TM)+w_2(E)\equiv c \bmod 2.
\end{equation}
\end{prop}
A proof of this proposition is in \cite{FT2000}.

\section{Proof of the Theorem}\label{s:proof}

Throughout the rest of the paper, \(M_1\) and \(M_2\) are closed
Riemannian manifolds. By virtue of the K\"unneth theorem \cite{Greenberg81},
every element \(\alpha\in\Hm^2(M_1\times M_2,\mathbb Z/2)\) admits
a decomposition
\begin{equation}\label{e:decomp}
\alpha=\alpha_{20}+\alpha_{11}+\alpha_{02},\quad\alpha_{ij}\in\Hm^i(M_1,
\mathbb Z/2)\otimes\Hm^j(M_2,\mathbb Z/2).
\end{equation}
\begin{lemma}
An element \(\alpha\in\Hm^2(M_1\times M_2,\mathbb Z/2)\), decomposed as
in \eqref{e:decomp}, is the
\(\bmod\;2\) reduction of \(c\in\Hm^2(M_1\times M_2,\mathbb{Z})\) if,
and only if, there are elements \(c_{ij}\in\Hm^i(M_1,\mathbb{Z})
\otimes\Hm^j(M_2,\mathbb{Z})\) such that \(\alpha_{ij}\equiv c_{ij}\bmod 2\).
\end{lemma}
\begin{proof}
From the second sequence in Section \ref{ss:UCT} and the fact that
\(\Hm_0(M,\mathbb{Z})\) is free there follows
\begin{equation}\label{e:H1Hom}
\Hm^1(M,R)\cong\Hom(\Hm_1(M,\mathbb{Z}),R).
\end{equation}
Therefore,  the
groups \(\Hm^i(M,\mathbb{Z})\) are free for \(i=0\) and \(1\); the
K\"unneth theorem applied to \(c\) gives
\[
c=c_{20}+c_{11}+c_{02},\quad c_{ij}\in\Hm^i(M_1,\mathbb{Z})\otimes\Hm^j(M_2,
\mathbb{Z}).
\]
To complete the proof, one reduces both sides of the last equality modulo 2 and
notes that such a reduction respects the decomposition and acts on each summand
separately.
\end{proof}
The Whitney formula gives
\[
\;w_2(T(M_1\times M_2))=w_2(TM_1)\otimes 1+w_1(TM_1)\otimes w_1(TM_2)
+1\otimes w_2(TM_2).
\]
Applying Lemma 1 to the last equation one obtains that a {\it necessary\/}
condition for the existence of a pin\(^c\) structure on \(M_1\times M_2\)
is that both \(M_1\) and \(M_2\) have pin\(^c\) structures. The condition
becomes {\it sufficient\/} if, in addition, there holds:
\begin{eqnarray}\label{e:cond}
\alpha_{11}=
w_1(TM_1)\otimes w_1(TM_2)\quad\mbox{\it is the {\rm mod 2} reduction}\\
\mbox{\it of an element}\quad
 c_{11}\in\Hm^1(M_1,\mathbb{Z})\otimes\Hm^1(M_2,\mathbb{Z}).\nonumber
\end{eqnarray}
To complete the {\it proof of the Theorem}, one has to show that
the alternative `(i) or (ii)' is equivalent to \eqref{e:cond}.
Part (i) is easy: one of the manifolds under consideration is orientable, if and
only if, one can take \(c_{11}=0\) in \eqref{e:cond}. Assume now
that (i) does not hold so that \(\alpha_{11}\neq 0\). Since our
manifolds are closed, their homology and cohomology groups are
finitely generated,
\begin{eqnarray}
\Hm_1(M_1,\mathbb{Z})&\cong \mathbb{Z}^k\oplus\mathbb{Z}/2^{p_1}\oplus
\dots\oplus\mathbb{Z}/2^{p_s}\oplus T_1,\\ \label{e:1}
\Hm_1(M_2,\mathbb{Z})&\cong \mathbb{Z}^l\oplus\mathbb{Z}/2^{q_1}\oplus
\dots\oplus\mathbb{Z}/2^{q_t}\oplus T_2,\label{e:2}
\end{eqnarray}
where \(T_i\) is the torsion of \(\Hm_1(M_i,\mathbb{Z})\) other than
the 2-torsion part so that \(T_i\otimes_{\mathbb{Z}}\mathbb Z/2=0\).
Using \eqref{e:H1Hom}
one can write
\begin{eqnarray*}
\Hm^1(M_1,\mathbb{Z}) =&  \mathbb{Z}^k\quad\qquad\;\,
& {\mbox{with basis}}\; (a_1,\dots,a_k),\\
\Hm^1(M_2,\mathbb{Z}) =& \mathbb{Z}^l\quad\qquad\;\,
&{\mbox{with basis}}\; (b_1,\dots,b_l),\\
\Hm^1(M_1,\mathbb Z/2) =& \mathbb Z/2^k\oplus\mathbb Z/2^s\quad
&{\mbox{with basis}}\; (a'_1,\dots,a'_k,c'_1,\dots,c'_s),\\
\Hm^1(M_2,\mathbb Z/2) =& \mathbb Z/2^l\oplus\mathbb Z/2^t\quad
&{\mbox{with basis}}\; (b'_1,\dots,b'_l,d'_1,\dots,d'_t),
\end{eqnarray*}
where \(a'_i\) and \(b'_i\) is the reduction mod 2 of \(a_i\)
 and \(b_i\), respectively. The classes \(w_1(TM_1)\)
and \(w_1(TM_2)\) can be decomposed with respect to the above bases
and one sees that condition \eqref{e:cond} is satisfied if, and only if,
\(w_1(TM_1)\) and \(w_1(TM_2)\) do not contain summands with
\(c'_i\) and \(d'_i\), respectively.\qed
\section{Examples}\label{s:ex}
\subsection{Products of pin\(^c\) manifolds}
\begin{exmp}\label{ex:1}
Recall that the real projective space \(\RP_{2k}\), \(k\in\mathbb{N}\),
has a pin structure \cite{DT86}; therefore, {\it a fortiori}, it has
a pin\(^c\) structure.
If   \(M_1=\RP_{2k}\) and \(M_2=\RP_{2l}\),
where \(k\) and \(l\) are positive integers, then
 \(M_1\times M_2\) has no pin\(^c\) structure.
\end{exmp}
\begin{proof}
Indeed, \(w_1(TM_1)\otimes w_1(TM_2)\neq 0\)
cannot be the reduction of an integral cohomology element because
\(\Hm^1(\RP_n,\mathbb{Z})=0\) for \(n\geqslant 2\).
\end{proof}
\begin{exmp}
For every integer \(k\geqslant 5\) we now construct a
non-orientable manifold
\(M_k\) of dimension \(k\) so that \(M_k\times M_k\) has a pin\(^c\)
structure. Let \(B\) be the two-dimensional M\"obius band
and let \(D_k\) be the unit ball of dimension \(k\) so that its boundary is
\(\mathbb{S}_{k-1}\). The manifold \(M_k\) is defined as the sum of two
manifolds, glued along their common boundary,
\[
M_k=(B\times\mathbb{S}_{k-2})\cup_{\mathbb{S}_1\times\mathbb{S}_{k-2}}
(\mathbb{S}_1\times D_{k-1}).
\]
The non-orientability of \(B\) implies that of \(M_k\). One easily finds
\[
\Hm_0(M_k,\mathbb{Z})=\mathbb{Z},\quad\Hm_1(M_k,\mathbb{Z})=\mathbb{Z},\quad
\Hm_2(M_k,\mathbb{Z})=0\quad\mbox{for \(k\geqslant 5\)}.
\]
These groups are free and so  are the groups
\(\Hm_i(M_k\times M_k,\mathbb{Z})\) for \(i=0,1,2\). This implies that
all elements of \(\Hm^2(M_k\times M_k,\mathbb{Z}/2)\) have an integral
lift; therefore, the manifold \(M_k\times M_k\) has a pin\(^c\) structure, but,
being a product of two non-orientable manifolds, it has no pin structure.
\end{exmp}
\subsection{Lipschitz manifolds without pin\(^c\) structure}
In this section, we construct a series of examples of manifolds with
a Lipschitz structure and without pin\(^c\) structure.
Following a remark in \cite{KT90}, we show that
every three-dimensional manifold is pin\(^-\). We  construct examples
of five-dimensional manifolds with a Lipschitz structure
 and without pin\(^c\) structure.
\begin{exmp}\label{ex:lnonpinc}
From Example \ref{ex:1}, it follows that \(M=\RP_2\times\RP_2\)
has no pin\(^c\) structure. Therefore, the same is true of the
five-dimensional manifold \(N=M\times\mathbb{S}_1\). We now construct
a vector bundle \(E\to N\), with fibres of real dimension 2, such
that \(TN\oplus E\) is a trivial. The embedding \(\RP_2\to\RP_3\cong\SO_3\)
induces on \(\RP_2\) a trivial bundle so that there is a line bundle
\(l\) such that \(T\RP_2\oplus l\) is trivial. One puts
\(E=l\oplus l\).  From the Whitney product
theorem one obtains \(w_1(E)=w_1(TN)\) and
\[
w_2(TN)+w_2(E)=w_1(TN)^2.
\]
Since, for every manifold \(M\),  \(w_1(TM)^2\) is the reduction mod~2 of an
integral element (see, e.g., Section 6 in \cite{FT2000}), the criterion
\eqref{e:FT} for the existence of a Lipschitz structure on \(N\) is
satisfied.
\end{exmp}

If \(M\) is a {\it three\/}-dimensional manifold, then Wu's formula (see \S 11
in \cite{MiSt74}) gives
\begin{eqnarray}
w_1(TM) &=& \Sq^0(v_1) + \Sq^1(v_0) = v_1 \label{wuw1}\\
w_2(TM) &=& \Sq^0(v_2) + \Sq^1(v_1) + \Sq^2(v_0) = v_2 + (v_1)^2 \label{wuw2}\\
v_2 &=& 0. \label{wuv2}
\end{eqnarray}
Here \(v_i\in\Hm^i(M,\mathbb{Z}/2)\) is the \(i\)th Wu class of
\(M\) characterized by
\[
v_i\cup w=\Sq^i(w)\;\;\mbox{for every
\(w\in\Hm^{3-i}(M,\mathbb{Z}/2)\)}
\]
so that \(v_i=0\)
for \(i>3-i\) implying \eqref{wuv2}. Equations \eqref{wuw1}-\eqref{wuv2} give
now
 \(w_2(TM)=w_1(TM)^2\); this shows that
every three-dimensional manifold has pin\(^-\) structure;
therefore, {\it a fortiori\/}, a pin\(^c\) structure (Section \ref{ss:exist}). As a
consequence, five is the smallest dimension of a manifold admitting
a Lipschitz structure, but no pin\(^c\) structure.

Example \ref{ex:lnonpinc} can be generalized; it suffices to note
that the example is based on two facts: \(\RP_3\) is parallelizable
and \(T\mathbb{S}_1\) is stably trivial (even trivial).
The projective space \(\RP_7\) is also parallelizable
and there are many manifolds with a stable trivial tangent bundle (Lie groups,
spheres, orientable three-dimensional manifolds and products of these
manifolds). This leads to the following example:
\begin{exmp}
 If \(M\) is \((2n+1)\)-dimensional and such
that
\[
TM\oplus\theta^{2(k+l+1)}\cong\theta^{2(k+l+n+1)+1}\quad\mbox{where}
\quad k,l\in\{1,3\},
\]
then the manifold
\begin{equation}\label{e:manif}
\RP_{2k}\times\RP_{2l}\times M
\end{equation}
has a Lipschitz structure, but no pin\(^c\) structure. Indeed,
denoting by \(l_i\) the normal line bundle for the embedding
\(\RP_{2i}\to\RP_{2i+1}\), one can take
\[
E = l_k \oplus l_l.
\]
One easily checks that
\[
T(\RP_{2k}\times \RP_{2l}\times M) \oplus l_k\oplus l_l
\cong TM\oplus \theta^{2(k+l+1)}\cong \theta^{2(k+l+m+1)+1}.
\]
and so the manifold \eqref{e:manif} has the desired property.
\end{exmp}
\begin{acknowledgements}
One of us (A.T.) is grateful to Boris Dubrovin for a discussion in Trieste on
how to solve the problem of existence of pin\(^c\) structures on products
of manifolds. This research was supported in part by
the Polish Committee  for Scientific Research (KBN) under grant
2 P03B 060 17.
\end{acknowledgements}

\end{article}
\end{document}